\documentclass[11pt]{article}
\usepackage{mathptmx}
\usepackage{graphicx, amssymb, amsthm}
\usepackage[papersize={5.5in, 8.5in}, top=0.35in, bottom=0.5in, left=0.35in, right=0.35in]{geometry}
\usepackage[colorinlistoftodos]{todonotes}
\usepackage{float}
\usepackage{amsmath}
\usepackage{algpseudocode}
\usepackage{algorithmicx}
\usepackage{algorithm}

\usepackage{color}
\usepackage{tikz}
\usepackage{url}
\usepackage{hyperref}

\date{}

\theoremstyle{definition}
\newtheorem{df}{Definition} [section]

\theoremstyle{plain}
\newtheorem{thm}[df]{Theorem}

\newtheorem{problem}[df]{Problem}
\newtheorem{claim}[df]{Claim}

\title{A $6$-chromatic two-distance graph in the plane}
\author{%
Geoffrey Exoo\\
Department of Mathematics and Computer Science\\
Indiana State University\\
Terre Haute, IN 47809,
\texttt{ge@cs.indstate.edu}
\and
Dan Ismailescu\\
Mathematics Department\\
Hofstra University\\
Hempstead, NY 11549,
\texttt{dan.p.ismailescu@hofstra.edu}
}

\begin{document}

\maketitle
\thispagestyle{empty}
\pagestyle{empty}
\begin{abstract}
\noindent We prove that if one colors each point of the Euclidean plane with one of five colors, then there exist two points of the same color that are either distance $1$ or distance $2$ apart.
\end{abstract}


\maketitle
\section{\bf Introduction}

The classical Hadwiger-Nelson problem asks for the minimum number of colors needed to color the points of the Euclidean plane in such a way that no two points distance $1$ apart are identically colored.

This quantity is denoted by $\chi(\{1\})$ or simply $\chi$, and it is known as \emph{the chromatic number of the plane}. According to Soifer \cite{soifer}, the problem of determining $\chi$ was posed in the 1950's by Edward Nelson. It is easy to prove that $4\le \chi \le 7$. The lower bound follows from an example due to the Moser brothers who constructed a unit-distance graph with chromatic number $4$. An upper bound of $7$ can be proved by considering a $7$-coloring of the regular hexagon tiling of the plane with hexagons with diameter slightly less than $1$.

Rather surprisingly, no progress has been achieved until 2018 when de Grey \cite{degrey} constructed a unit-distance graph with chromatic number $5$. A different construction was obtained independently and at about the same time by the authors of this note \cite{exooismailescuchi5}.

De Grey's initial graph has $1581$ vertices. This was soon after improved by Heule who managed to reduce the order to $553$ vertices \cite{heule1} and subsequently to $529$ vertices \cite{heule2} . As of the time of this writing, the record is held by Parts \cite{parts} who found a $5$-chromatic unit-distance graph with $510$ vertices and $2508$ edges.

\noindent In this paper we consider the \emph{two forbidden distances} variation of the Hadwiger-Nelson question.
\begin{problem}
For a given number $d>1$, what is the minimal number of colors needed to color the plane such that no two identically colored points are either distance $1$ or distance $d$ apart?
\end{problem}

We denote this number by $\chi(\{1,d\})$.

The problem has been considered before; the earliest reference we could find is due to Owings, Tetiva and Huddleston \cite{HOT}. In this paper it is proved that $\chi(\{1,d\})\ge 5$ for $d\in \{\sqrt{2},\sqrt{3},(\sqrt{6}+\sqrt{2})/2, (\sqrt{5}+1)/2\}$.

A few years later, Katz, Krebs and Shaheen \cite{KKS}, unaware of the results of \cite{HOT}, gave a different proof of $\chi(\{1,\sqrt{2}\})\ge 5$. The list of values $d$ for which $\chi(\{1,d\})\ge 5$ was further extended by the present authors \cite{exooismailescu2distances}.

Of course, since $\chi(\{1,d\})\ge \chi(\{1\})$ for all values of $d$, all the above results are immediate consequences of de Grey's bound.

As far as we know, the only instance in which one can prove a better lower bound than $5$ is when $d=(\sqrt{5}+1)/2$. In this case Huddleston showed that
\begin{thm}\cite{HOT}\label{huddleston}
\begin{equation*}
\chi\left(\left\{1, \frac{\sqrt{5}+1}{2}\right\}\right)\ge 6.
\end{equation*}
\end{thm}

The following concept extends the notion of unit-distance graph to the case of two distances.
\begin{df}
A $\{1,d\}-graph$ is a graph whose vertices are points in the plane, with two vertices being adjacent if the Euclidean distance between them is either $1$ or $d$.
\end{df}

Huddleston's proof relies on the fact that a regular pentagon with side length $1$ has diagonals of length $(\sqrt{5}+1)/2$, that is, the complete graph $K_5$ can be represented as a $\{1,(\sqrt{5}+1)/2\}$-graph. It follows from results of Einhorn and Schoenberg \cite{ES} that this is the only case in which $K_5$ is a $\{1,d\}$-graph.

In this note we prove a similar result to that in Theorem \ref{huddleston} in the case $d=2$.

\begin{thm}\label{mainthm}
\begin{equation*}
\chi\left(\left\{1, 2\right\}\right)\ge 6.
\end{equation*}
\end{thm}

\section{\bf Proof of Theorem \ref{mainthm}}

Our approach to proving this result is constructive - we will present a finite $\{1,2\}$-graph which cannot be $5$-colored.
We build this graph in several stages.
As in the proofs of $\chi(\{1\})\ge 5$ in \cite{exooismailescuchi5, degrey}, we will use vertices with coordinates of the form
$(a\sqrt{3}/12+b\sqrt{11}/12, c/12+d\sqrt{33}/12)$ where $a, b, c$, and $d$ are integers.

In order to improve readability, we will use the following notation
\begin{equation}\label{notation}
[a,b,c,d]:=\left(\frac{a\sqrt{3}}{12}+\frac{b\sqrt{11}}{12}, \frac{c}{12}+\frac{d\sqrt{33}}{12}\right).
\end{equation}

{\bf Step 1.} Consider the following set of $23$ points in the plane:
\begin{align*}
S:=\{&[0, 0, 0, 0], [0, 0, 0, -4], [0, 0, -6, -2], [0, 0, -6, 2], [-6, 0, 0, -2], \\
     &[-4, 0, 0, 0],[-4, 0, -6, -2],[-4, 0, -6, 2], [-2, 0, 0, -2], [-2, 0, -6, -4],\\
     & [-2, 0, -6, 4],[0, -6, -6, 0],[-5, -3, 3, 3], [-5, 3, -3, 3], [-2, -6, 0, 0],\\
     & [-2, -6, 0, -4], [-2, 6, 0, 0], [-2, 6, 0, -4],[-6, -6, 0, 0], [-6, 6, 0, 0],\\
     & [-4, 0, 0, -4], [0, 0, -12, 0], [-8, 0, 0, 0]\}.
\end{align*}

Next, consider all the points in $S$ together with their reflections across the $x$-axis and across the $y$-axis, respectively.
One obtains a new set $T$ which has $57$ points. For $k=0..5$, let $U_k$ be the image of $T$ under a rotation of angle $k\pi/3$ about the origin, and define $V:=U_0 \cup U_1 \cup U_2 \cup U_3 \cup U_4 \cup U_5$.

It is easy to check that the image of $[a,b,c,d]$ under a rotation of $\pi/3$ about $[0,0,0,0]$ is $[(a-c)/2, (b-3d)/2,(3a+c)/2,(b+d)/2]$.
Let $G$ be the $\{1,2\}$-graph whose vertex set are the points in $V$.
\begin{claim}
The graph $G$ has $205$ vertices, $966$ edges of length $1$, $423$ edges of length $2$, and exactly $18$ $5$-colorings.
\end{claim}

{\bf Step 2.}
We construct a slightly larger graph $H$, by including the following nine additional vertices to the vertex set of $G$:
\begin{align*}
A:=&{\bf[-2, 0, 0, -6]}, B:={\bf [8, 0, 0, 4]}, [-4, -6, -6, -4], [-4, 6, 6, -4],\\
& [-3, -3, -3, -5], [-4, 0, -12, 4], [-4, 0, 12, 4], [7, -3, 3, 3], [7, 3, -3, 3].
\end{align*}

\begin{claim}
The graph $H$ defined above has $214$ vertices, $1004$ edges of length $1$, $446$ edges of length $2$, and exactly $35$ $5$-colorings. Moreover, in each of these colorings, vertices $A$ and $B$ are of the same color.
\end{claim}

\begin{figure}[H]
\centering
\includegraphics[width=0.97\linewidth]{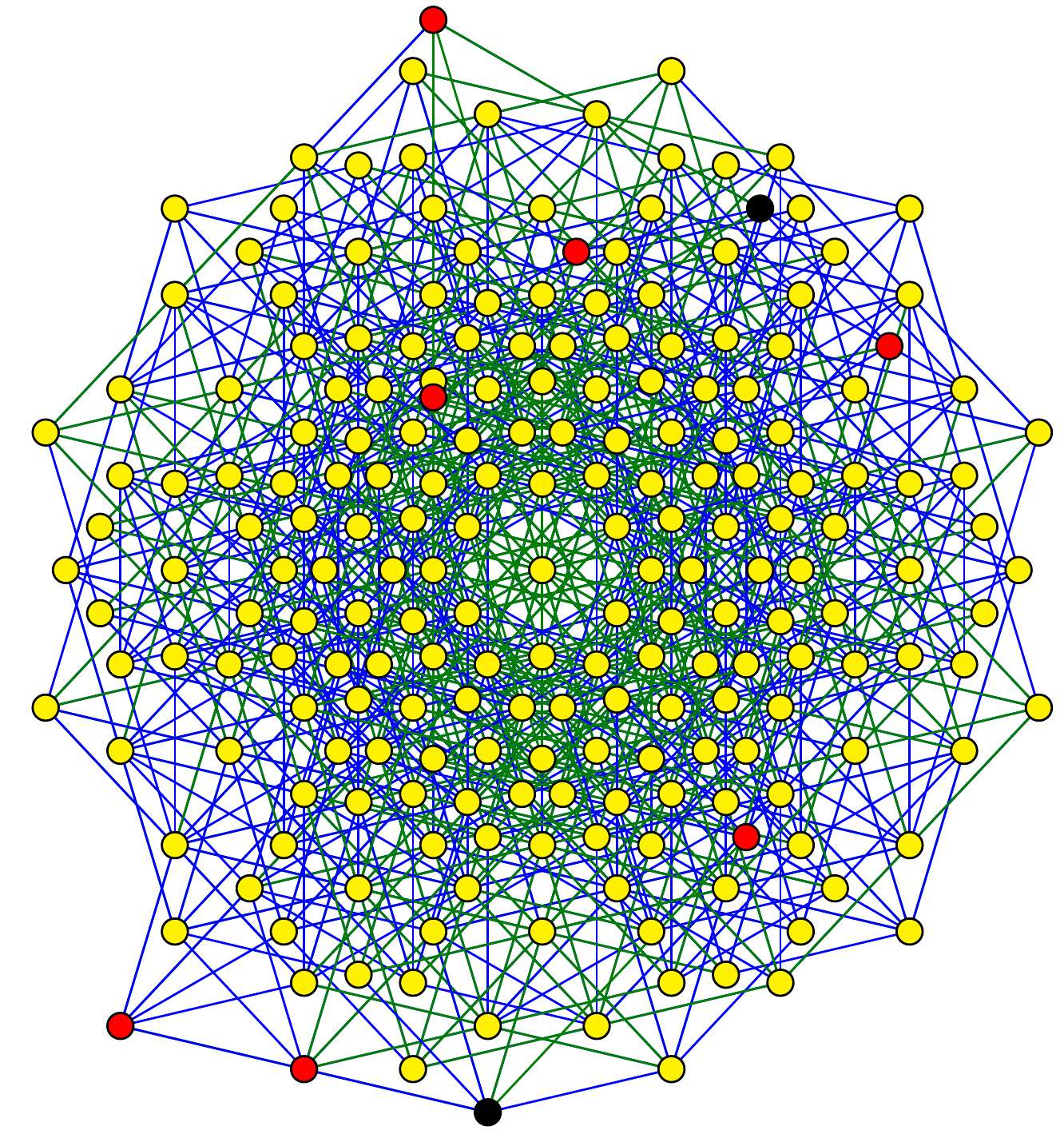}
\caption{The $\{1,2\}$-graph $H$. Vertices $A$ and $B$ appear in black, the other seven new vertices in red.}
\label{G214}
\end{figure}

We describe the techniques used to substantiate these claims in section 3.

Note that the distance between $A$ and $B$ is exactly $5$.
Rotating the vertices of $H$ about vertex $A$ by an angle $\arccos(49/50)=\arcsin(3\sqrt{11}/50)$, creates a copy of $H$, which we denote $H'$.
The image of vertex $B$ under this rotation is a point $B'\in V(H')$, and the distance between $B$ and $B'$ is exactly $1$.

Let $K$ be the $\{1,2\}$-graph whose vertex set is $V(H)\cup V(H')$.
It can be checked that $K$ has $426$ vertices, $2009$ edges of length $1$, and $892$ edges of length $2$. More importantly, every $5$-coloring of $K$ forces vertices $A$, $B$, and $B'$ to receive the same color. Since $AB=AB'=5$ and $BB'=1$, it follows that $\chi(K)\ge 6$. This concludes the proof of Theorem \ref{mainthm}.

It is interesting to note that the final argument in the proof of Theorem \ref{huddleston} involves an isosceles triangle of side lengths $5, 5, 1$, exactly the same triangle as in our construction.

All vertices of the graph $K$ have coordinates in $\mathbb{Q}[\sqrt{3},\sqrt{11}]$. This is because the smallest angle of an isosceles triangle with sides $5, 5, 1$ is $\arccos(49/50)=\arcsin(3\sqrt{11}/50)$.

\section{\bf Computations}\label{comp}

In this section we describe the methods used to establish the coloring counts
claimed above.
The assertions pertaining to edge counts for both of the graphs can be easily obtained
by direct computation using the data available at \cite{webvertices,webmatrix}.
The assertions that there are exactly $18$ $5$-colorings for graph $G$,
and $35$ $5$-colorings for graph $H$ require more difficult computations.

To find all $5$-colorings for these graphs we used a simple recursive exhaustive search procedure
that allowed us to divide the work across multiple processors.  The outline of the search
procedure is given below.  Before
the procedure is used, the vertices are ordered as follows.

\begin{itemize}
\item The vertices are partitioned into orbits based on the dihedral group generated by
the transformations used in the construction (reflections in the axes and the $\pi/3$ rotation).
\item Vertices within an orbit are sorted by polar angle: $0 \leq \theta < 2 \pi$, and at all
stages appear consecutively in the vertex ordering.
\item Vertex orbits are sorted in descending order by degree.
\item In case of ties, vertex orbits adjacent to the largest number of vertices that appear earlier in
the ordering are listed first.
\end{itemize}

Then each vertex is assigned the $NC$ (uncolored) value, and the following search procedure is called
with vertex $0$ and the list of colors as parameters.

\begin{algorithm}[htp]

\caption{Coloring Search}
\label{coloralgo}
\begin{algorithmic}
\algnewcommand{\BigComment}[1]{\State \(\triangleright\) #1}
\Procedure{Search}{$vertex,colors$} \Comment{The graph is global}
  \State Declare static variable $calls = 0$

  \If{$vertex = n$} \Comment{Vertices numbered $0 \cdots n-1$}
    \State Output $colors$
    \State Return
  \EndIf

  \State
  \BigComment{Next we split the work among $NCPU$ processors.  Each processor has an id ($cpu$).
           The splitting is done at level $keydepth$ (where $15 \leq keydepth \leq 20$ here).}
  \State

  \If{$vertex = keydepth$} \Comment{Level where splitting is done}
    \State $calls \gets calls + 1$
    \If {$calls \not\equiv cpu \pmod{NCPU}$}
      \State Return
    \EndIf
  \EndIf

  \State
  \BigComment{If the color of a vertex has already been forced, we skip ahead to the next level.}
  \State

  \If{$color[vertex] \not= NC$} \Comment{Vertex already colored}
    \State Search($vertex+1,colors$)
    \State Return
  \EndIf
  \State
  \BigComment{The $Assign$ function colors a vertex and recursively considers all implications.
    This may force colors on other vertices, or conflicts may be discovered, in which case it returns False.
    The $UnAssign$ function restores the previous state of the coloring.}
  \State
  \For{each color $c$}
    \If{Assign($vertex, colors, c$)}
      \State Search($vertex+1,colors$)
      \State UnAssign($vertex, colors, c$)
    \EndIf
  \EndFor
\EndProcedure
\end{algorithmic}
\end{algorithm}

The computations were performed using 48 threads on an
{\em AMD EPYC 7551 32-Core (64 Virtual Core) Processor},
and were completed in $3780$ seconds of elapsed time and $81000$ seconds of total processing time for $G$,
and $5120$ seconds of elapsed time and $95000$ seconds of total processing time for $H$.
In each case, all but three of the threads were finished halfway through the computation, which was not
surprising, given that our method for splitting the work was fairly crude.

\end{document}